\documentclass[twocolumn,aps,secnumarabic,nofootinbib,superscriptaddress,%
    floatfix]{revtex4}
\usepackage{amssymb,amsmath2000,amsthm2000,times,mathptmx}
\newtheorem{theorem}{Theorem}
\newtheorem{corollary}{Corollary}
\newtheorem{algorithm}{Algorithm}
\newtheorem*{claim}{Claim}
\newtheorem*{conjecture}{Conjecture}
\theoremstyle{remark}
\newtheorem{example}{Example}

\newcommand{\comment}[1]{} % comment out previous line for final version
\newcommand{\muspace}[1]{\mskip #1mu plus 1mu minus 1mu}
\newcommand{\deq}{\muspace{5} = \muspace{5}}
\newcommand{\drel}[1]{\muspace{5} #1 \muspace{5}}

\newcommand{\opshrink}[2]{\mskip -#1mu minus -#1mu #2 \mskip -#1mu minus -#1mu}
\newcommand{\squish}[1]{\opshrink{3}{#1}}
\newcommand{\bcdisplaydiv}{\muspace{-4} \Bigg/ \muspace{-4}}
\newcommand{\textbc}[2]{\mbox{\small $\binom{#1}{#2}$}}

\newcommand{\AG}{\mathop{\rm AG}}
\newcommand{\GF}{\mathop{\rm GF}}

\newenvironment{algbody}{%
        \vspace{-3.75ex plus 1ex minus -1ex}\nopagebreak
        \begin{enumerate}%
        \setlength{\topsep}{0ex}%
        \setlength{\itemsep}{0ex plus .5ex minus .5ex}%
        }{\end{enumerate}%
          \vspace{-1ex plus .5ex minus .5ex}}
\newenvironment{subalgbody}{%
        \vspace{-1.5ex plus .5ex minus .5ex}
        \begin{enumerate}%
        \setlength{\topsep}{0ex}%
        \setlength{\itemsep}{0ex plus .5ex minus .5ex}%
        }{\end{enumerate}}

\begin{document}

\title{Asymptotically optimal covering designs}
\author{Daniel M. Gordon}
\affiliation{Center for Communications Research,
    4320 Westerra Ct., San Diego, CA 92121}
\author{Greg Kuperberg}
\affiliation{Department of Mathematics,
    Yale University, New Haven, CT 06520}
\author{Oren Patashnik}
\affiliation{Center for Communications Research,
    4320 Westerra Ct., San Diego, CA 92121}
\author{Joel H. Spencer}
\affiliation{Courant Institute, NYU, New York, NY 10012}

% \date{10 November 1995}

\begin{abstract}
A $(v,k,t)$ \emph{covering design}, or \emph{covering}, is a family of
$k$-subsets, called \emph{blocks}, chosen from a $v$-set, such that each
$t$-subset is contained in at least one of the blocks. The number of blocks is
the covering's \emph{size}, and the minimum size of such a covering is denoted
by $C(v,k,t)$. It is easy to see that a covering must contain at least
$\binom{v}{t}/\binom{k}{t}$ blocks, and in 1985 R\"odl~\cite{rodl} proved a
long-standing conjecture of Erd\H{o}s and Hanani~\cite{erdos-hanani} that for
fixed $k$~and~$t$, coverings of size $\binom{v}{t}/\binom{k}{t} (1+o(1))$
exist (as $v \to \infty$).

An earlier paper by the first three authors~\cite{gkp} gave new methods for
constructing good coverings, and gave tables of upper bounds on $C(v,k,t)$ for
small $v$,~$k$, and~$t$. The present paper shows that two of those
constructions are asymptotically optimal: For fixed $k$~and~$t$, the size of
the coverings constructed matches R\"odl's bound. The paper also makes the
$o(1)$ error bound explicit, and gives some evidence for a much stronger bound.
\end{abstract}

\maketitle

\section{Introduction}

Let the covering number $C(v,k,t)$ denote the smallest number of
$k$-subsets of a $v$-set that cover all $t$-subsets.
The best general lower bound on $C(v,k,t)$,
due to Sch\"onheim~\cite{schonheim}, comes from the following inequality:
\begin{theorem}
\label{schonheim-one-level}
$$C(v,k,t) \drel{\geq} \Bigl\lceil \frac{v}{k} \,
C(v\squish{-}1, k\squish{-}1, t\squish{-}1) \Bigr\rceil \,.$$
\end{theorem}
Iterating this gives the Sch\"onheim bound
$$ C(v,k,t) \drel{\geq} \Bigl\lceil
\frac{v}{k} \Bigl\lceil \frac{v-1}{k-1} \ldots
\Bigl\lceil \frac{v-t+1}{k-t+1} \Bigl\rceil \ldots \Bigl\rceil \Bigl\rceil \,.
$$
The best general upper bound on $C(v,k,t)$ is due to R\"odl~\cite{rodl}:
The \emph{density\/} of a covering is the
average number of blocks containing a $t$-set.
The minimum density is $C(v,k,t) \textbc{k}{t} / \textbc{v}{t}$, and
is obviously at least~1.  R\"odl showed that for $k$~and~$t$ fixed
there exist coverings with density $1+o(1)$ as $v$~gets large.

This paper shows that two of our constructions~\cite{gkp}
match the bound of R\"odl's theorem.
One of the constructions gives an easier proof of the theorem
than R\"odl's original proof \cite{rodl}.
The other construction provides a computationally efficient version
of R\"odl's theorem.
In Section~\ref{s:review} we review the two constructions.
In Section~\ref{s:greedy} we show that the first one,
which uses a greedy algorithm, is asymptotically optimal.
And in Section~\ref{s:induced} we show that the second one,
which constructs an induced covering from a finite-geometry covering,
is also asymptotically optimal,
and that it is computationally efficient as well.

Theorem~\ref{t:greedy} (in Section~\ref{s:greedy})
is a special case of a main result of the fourth author~\cite{spencer};
R\"{o}dl and Thoma~\cite{rodl-thoma} gave another proof of that result.
We present the proof here to keep the paper self-contained
and to provide an explicit error bound for use in Section~\ref{s:induced}.

\section{Covering Constructions}
\label{s:review}

Here we summarize two methods
for constructing asymptotically optimal coverings.
Our previous paper~\cite{gkp} gives more details,
as well as computational results for small $v$,~$k$, and~$t$.

\subsection{Greedy Coverings}
\label{ss:greedy}

\begin{algorithm}
\label{a:greedy}
{\rm Random Greedy $(v,k,t)$ Covering}
\end{algorithm}
\begin{algbody}

\item
\label{ordering}
Fix a random ordering of the $k$-sets of a $v$-set.

\item
\label{main-step}
Choose the earliest $k$-set containing no already-covered $t$-set.

\item
\label{repeat}
Repeat Step~\ref{main-step} until no $k$-set can be chosen.

\item
\label{the-coverup}
Cover the remaining $t$-sets with one $k$-set each.
\end{algbody}

This greedy algorithm is a little different from our previous one. That
algorithm uses one of four possible orderings in Step~\ref{ordering}:
lexicographic, colex, Gray code, or random. Also, it chooses in
Step~\ref{main-step} the earliest $k$-set that contains the most
still-uncovered $t$-sets; thus it continues with Steps
\ref{main-step}~and~\ref{repeat} instead of cutting out to
Step~\ref{the-coverup}. That algorithm produces slightly better coverings in
practice, but is harder to analyze than the algorithm here.

\subsection{Induced Finite Geometry Coverings}
\label{s:geometries}

The $k$-flats of an affine or projective geometry form a covering. For this
paper, we restrict our attention to the hyperplanes of an affine geometry,
which form an optimal covering:
\begin{theorem}
For a prime power~$q$ and integer $t>1$,
the hyperplanes of the affine geometry $\AG(t,q)$
are a $(q^t\!, q^{t-1}\!, t)$ covering of size
$$
C(q^t\!, q^{t-1}\!, t) \deq \frac{q^{t+1}-q}{q-1} \,.
$$
\end{theorem}
The density of such a covering is
$$
\frac{q^{t+1}-q}{q-1} \binom{q^{t-1}}{t} \bcdisplaydiv \binom{q^t}{t}
     \deq 1 + O(q^{-1}) \,.
$$

\begin{algorithm}
\label{a:induced}
{\rm Induced $(v,k,t)$ Covering}
\end{algorithm}
\begin{algbody}

\item
\label{induced-choose}
Choose a prime~$p$ with $p^t>v$,
and an integer~$\ell$, as specified later.

\item
\label{induced-precompute}
Precompute $(\ell',k,t)$ coverings, for $\ell < \ell'< 9 \ell$,
using Algorithm~\ref{a:greedy}.

\item
Choose $v$~points of the $\AG(t,p)$ at random.

\item
For each hyperplane, find its intersection with the $v$~points;
let $\ell'$ be the size of the intersection.
\begin{subalgbody}

\item
\label{main-case}
If $\ell <  \ell' < 9 \ell$,
add the blocks of the $(\ell',k,t)$ covering on those points
to the $(v,k,t)$ covering.

\item
If $\ell' \leq \ell$ or $\ell' \geq 9 \ell$,
trivially add $\textbc{\ell'}{t}$ blocks to the $(v,k,t)$ covering.

\end{subalgbody}
\end{algbody}

The new blocks each have $k$~elements, and together they cover all $t$-sets,
so they form a $(v,k,t)$ covering.
The blocks of the affine covering and their intersection with the
$v$-set may quickly be computed by solving linear equations over $\GF(q)$.

This construction, too, differs slightly from our earlier version~\cite{gkp}.
In that paper, we construct $(\ell',k,t)$ coverings for all $\ell' < v$
by whatever construction gives the best results,
and then always use Step~\ref{main-case}.
That results in better coverings in practice, but is harder to analyze.

\section{Greedy Coverings and R\"odl's Bound}
\label{s:greedy}

The usual proofs of R\"odl's theorem
(R\"odl~\cite{rodl} or Alon and Spencer~\cite{alon-spencer})
seem nonconstructive;
however, they are actually analyses of a covering algorithm,
similar to the greedy algorithm with random ordering,
that constructs a covering in two steps.
First, it chooses a sequence of R\"odl nibbles,
each of which is a small, random collection of $k$-sets
that do not contain any $t$-set contained in any previous R\"odl nibble.
Second, when there is no longer room for a nibble,
it chooses a separate $k$-set for each remaining uncovered $t$-set.

The main difference between
the $k$-sets chosen in the sequence of R\"odl nibbles
and those chosen by the greedy algorithm
in Steps \ref{main-step}~and~\ref{repeat}
is that two $k$-sets in the same nibble may intersect each other in a $t$-set.
This difference seems small, hence it is natural to
conjecture that the greedy algorithm, too, meets R\"odl's bound.
It does:
\begin{theorem}
For fixed $k$~and~$t$, the greedy algorithm with random ordering
produces a covering with expected density $1+o(1)$ as $v \to \infty$.
\label{t:greedy}
\end{theorem}

The proof of Theorem~\ref{t:greedy} will proceed in several steps,
along the lines of Spencer~\cite{spencer}.

\subsection{The Continuous Model}

Model the execution of the greedy algorithm as a Poisson process; that is, a
given $k$-set is chosen between time $\tau$ and $\tau + \delta$ with
probability asymptotic to $\delta/\textbc{v-t}{k-t}$ as $\delta \to 0$, and the
probabilities of any two $k$-sets being chosen in any two time intervals are
independent. The process begins at time~$0$ and lasts forever. If a $k$-set
chosen by the process at some time~$\tau$ contains any previously covered
$t$-set, it \emph{fails\/} at time~$\tau$, otherwise it \emph{succeeds\/} and
its $t$-sets are considered covered after time~$\tau$. The $k$-set thus fails
at any time subsequent to~$\tau$ it is chosen.

The ordering determined by the first-choosings of the $k$-sets in this process
corresponds to the random ordering of the $k$-sets in the greedy algorithm, and
the $k$-sets that have succeeded at time infinity correspond to the $k$-sets
chosen by the greedy algorithm just prior to Step~\ref{the-coverup}. Thus to
prove the theorem it suffices to show that, at time infinity of the Poisson
process, a given $t$-set is covered with probability asymptotic to~1. (Since if
the proportion of $t$-sets covered at that point of the greedy algorithm goes
to~1 then so does the density of the eventual covering.) We actually find the
limit of this probability as $v \to \infty$ for every fixed~$\tau$, and we show
that this limit goes to~1 as $\tau \to \infty$.

Fix a time~$\tau$ and a $t$-set~$T$. Based on the Poisson process above, we
either define the \emph{dependence tree\/} of $(\tau,T)$ or else declare it to
be \emph{aborted}. The tree is rooted, and has $t$-vertices and
$k$-vertices---% $(\tau',T')$ with $|T'|=t$ and $(\tau',K')$ with $|K'|=k$. We
begin at time~$\tau$ with the tree consisting only of its $t$-vertex root
$(\tau,T)$, and we examine $k$-sets chosen by the process, proceeding backwards
in time from~$\tau$ toward~0.

There are three cases for a $k$-set~$K^*$ chosen at some time~$\tau^*$: if
$K^*$ does not contain any $T'$ already in the tree then do nothing; if it
contains two or more such $T'$ then declare the tree to be aborted; if (the
important case) it contains precisely one such $T'$ then add $(\tau^*,K^*)$ as
a child of $(\tau',T')$ and for every $t$-set $T^*\subset K^*$ \emph{
except\/}~$T'$ add $(\tau^*,T^*)$ as a child of $(\tau^*,K^*)$. We will say
that $T$ has given birth to $K^*$ at time $\tau^*$, and $K^*$ immediately gives
birth to all the $T'$ nodes.

The tree, if defined, is finite; a child of a $t$-vertex is a $k$-vertex and
vice versa. We label each vertex as follows. A $t$-vertex is \emph{covered\/} if
at least one of its children is accepted, else it is \emph{uncovered}; a
$k$-vertex is \emph{accepted\/} if none of its children is covered, else it is
\emph{rejected}. Thus a childless (leaf) $t$-vertex is uncovered, and a unique
labeling is defined inductively from the leaves up.

\begin{figure*}
\begin{center}
\setlength{\unitlength}{1pt}
\begin{picture}(320,200)(-150,0)
\newcommand{\toval}{\oval(56,18)}
\newcommand{\koval}{\oval(72,18)}
\newcommand{\tnode}[2]{\makebox(0,0){\footnotesize $#1,\muspace{6}\{#2\}$}}
\newcommand{\knode}[2]{\makebox(0,0){\footnotesize $#1,\muspace{10}\{#2\}$}}
\newcommand{\noderightlabel}[1]{\makebox(0,0)[l]{%
                                        \raisebox{-.5ex}[0ex][0ex]{\small #1}}}
\newcommand{\nodeleftlabel}[1]{\makebox(0,0)[r]{%
                                        \raisebox{-.5ex}[0ex][0ex]{\small #1}}}
\put(0,200){\toval} \put(0,200){\tnode{4.3}{1,2}}
\put(34,200){\noderightlabel{uncovered}}
\put(0,150){\koval} \put(0,150){\knode{3.7}{1,2,3}}
\put(42,150){\noderightlabel{rejected}}
\put(0,191){\line(0,-1){32}}
\put(-75,100){\toval} \put(-75,100){\tnode{3.7}{1,3}}
\put(-109,100){\nodeleftlabel{uncovered}}
\put(-14.25,140.5){\line(-3,-2){47.25}} % fudge factor because LaTeX screws up
\put(75,100){\toval} \put(75,100){\tnode{3.7}{2,3}}
\put(109,100){\noderightlabel{covered}}
\put(14.25,140.5){\line(3,-2){47.25}}
\put(75,50){\koval} \put(75,50){\knode{1.2}{2,3,4}}
\put(117,50){\noderightlabel{accepted}}
\put(75,91){\line(0,-1){32}}
\put(0,0){\toval} \put(0,0){\tnode{1.2}{2,4}}
\put(-34,0){\nodeleftlabel{uncovered}}
\put(60.75,40.5){\line(-3,-2){47.25}}
\put(150,0){\toval} \put(150,0){\tnode{1.2}{3,4}}
\put(116,0){\nodeleftlabel{uncovered}}
\put(89.25,40.5){\line(3,-2){47.25}}
\end{picture}
\end{center}
\caption{Example of a dependence tree.}
\label{example}
\end{figure*}
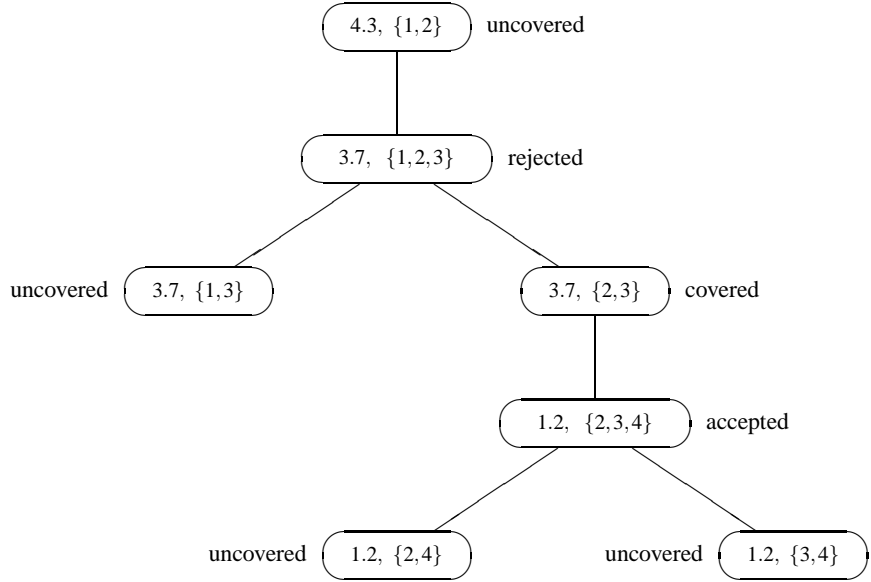

\begin{example}
Take $t\squish{=}2$; $k\squish{=}3$; $v\squish{=}10^{10}$; $\tau\squish{=}4.3$;
$T\squish{=}\{1,2\}$. Suppose $\{1,2,3\}$ is chosen at time~$3.7$ and
$\{2,3,4\}$ at time~$1.2$ and these are the only relevant chosen sets. The
dependence tree of $(4.3,\{1,2\})$ is shown in Figure~\ref{example}. Two of the
leaves $(1.2,\{2,4\})$ and $(1.2,\{3.4\})$ are uncovered, thus their parent
$(1.2,\{2,3,4\})$ is accepted, so $(3.7,\{2,3\})$ is covered and
$(3.7,\{1,2,3\})$ is rejected and finally $(4.3,\{1,2\})$ is uncovered. In the
corresponding Poisson process, $\{2,3,4\}$ succeeds at time~1.2, thus
$\{1,2,3\}$ fails at time~3.7, so no $3$-set covering $\{1,2\}$ is accepted by
time~$4.3$.
\end{example}

This example is consistent with the claim below.

\begin{claim}
Suppose the dependence tree of $(\tau,T)$ for some $\tau$~and~$T$ is defined.
Then $(\tau,T)$ is covered if and only if $T$~is covered by the Poisson
process.
\end{claim}

\begin{proof}[Proof of claim]
If $T$~is covered in the Poisson process by a $k$-set~$K$, then $K$~succeeded
at some time~$\tau^* \!$. Thus no $k$-set containing any of the $t$-sets
covered by~$K$ was chosen before~$\tau^* \!$, and $(\tau^*,K)$ is accepted,
hence $(\tau,T)$ is covered. Conversely, suppose that $(\tau,T)$ is covered in
its dependence tree. Then it has an accepted child. It might have several
accepted children, but since the tree is defined, the $k$-sets of these
children can intersect only in~$T$. The earliest such $k$-set succeeded, so it
covers~$T$. That establishes the claim.
\end{proof}

\subsection{The Idealized Tree}

The process above is still difficult to analyze directly, so we will define for
a fixed~$\tau$ an idealized process and an idealized tree, analogous to the
Poisson process and dependence tree. We will show that the idealized trees
behave like the dependence trees, and then find the probability that the root
of an idealized tree is covered.

The idealized tree has $t$-vertices and $k$-vertices, and consists at
time~$\tau$ just of a $t$-vertex root. Again, time goes backwards, from~$\tau$
to~0. In the interval from $\tau_1$ to $\tau_1-\delta$ each $t$-vertex has
probability asymptotic to~$\delta$ of giving birth to a $k$-vertex, which then
instantly gives birth to $D=\textbc{k}{t}-1$ new $t$-vertices. In a
length~$\delta$ interval each $t$-vertex has on average $\delta D$
grandchildren (also $t$-vertices), so the expected number of $t$-vertices goes
up by a factor of $1+\delta D$. The expected number of $t$-vertices at time~0
is thus $(1+\delta D)^{\tau/\delta} = e^{\tau D} (1 + O(\delta))$ as $\delta
\to 0$, hence with probability~1 the idealized tree is finite. The notions of
covered, uncovered, accepted, and rejected are defined on it as before.

We claim that the limit distribution of the dependence tree of $(\tau,T)$ as $v
\to \infty$ is the distribution for the idealized tree. Consider a fixed
idealized tree at time~$\tau$, and look at the dependence tree of $(\tau,T)$
from time $\tau_1$~to $\tau_1-\delta$ given that at~$\tau_1$ it matches the
idealized tree. The number of $t$-sets in the tree is $O(e^{\tau D})$, with
probability asymptotic to~1, so the number of $k$-sets that contain more than
one $t$-set already in the tree is $O(e^{2\tau D} v^{k-t-1})$, and thus the
probability of aborting (i.e., that some such $k$-set is chosen) is $O(\delta
e^{2\tau D} v^{-1})$. Therefore the total chance of aborting throughout the
length~$\tau$ interval is $O(\tau e^{2\tau D} v^{-1}) = o(1)$ for $\tau \leq
(\ln v)/(2+\epsilon) D$, for any fixed $\epsilon>0$.

For each $T'$ in the tree, the number of $k$-sets that contain $T'$ and no
other $t$-set in the tree is asymptotically $\textbc{v-t}{k-t}$, so $T'$~has a
($k$-vertex) child with probability asymptotic to~$\delta$, as in the idealized
version. Hence the two distributions are the same, as claimed.

Now we compute the probability $P(\tau)$ that the root of an idealized tree at
time~$\tau$ is uncovered. In the interval from $\tau$ to~$\tau-\delta$ of an
idealized process, a $t$-vertex either does or does not give birth, with
probabilities asymptotic to $\delta$ and $1-\delta$ as $\delta \to 0$. In the
former case, a $k$-vertex child is accepted with probability $P(\tau-\delta)^D
\!$, because each $t$-vertex grandchild has independent probability
$P(\tau-\delta)$ of being uncovered at time~$\tau-\delta$, and thus is rejected
with probability $1-P(\tau-\delta)^D \!$. Hence
\[
P(\tau) \drel{\sim}
\delta (1 - P(\tau-\delta)^D) P(\tau-\delta) + (1 - \delta) P(\tau-\delta) \,.
\]
So
$P(\tau-\delta) - P(\tau) \drel{\sim} \delta P(\tau-\delta)^{D+1} \!$,
which leads to the differential equation $P(\tau)' = -P(\tau)^{D+1}$
with the initial condition $P(0) = 1$.
The solution is
$$
P(\tau) \deq (\tau D + 1)^{-1/D} \,.
$$
In particular $\lim_{\tau \to \infty} P(\tau) = 0$,
so the root of an idealized tree at time infinity
is covered with probability asymptotic to~1.
Therefore, at time infinity of the Poisson process,
a given $t$-set is covered with probability asymptotic to~1,
and Theorem~\ref{t:greedy} is established.

\subsection{Estimating The Error Term}

The proof above shows that the greedy covering is optimal, but we have not
estimated the error term.  We conclude this section by giving a weak estimate,
along with some evidence for a stronger conjecture.

Consider the state of the algorithm at time $\tau = O(\log v)$. First, notice
that at this time of the Poisson process, the expected number of $k$-sets
chosen is $O(v^t \log v)$. Thus in the greedy algorithm it suffices to examine
just $O(v^t \log v)$ random $k$-sets before cutting out to
Step~\ref{the-coverup}. It takes only $O(v^t \log v)$ expected time and
$O(v^t)$ space to generate those $k$-sets (Brassard and Kannan~\cite{fly}), so
this early abort strategy dramatically speeds up the algorithm, at negligible
cost to the density of the covering:

\begin{corollary}
\label{cor:abort}
The early-abort greedy algorithm produces a covering with expected density
$1+o(1)$ in time $O(v^t \log v)$.
\end{corollary}

Second, at time $\tau = (\ln v)/(2+\epsilon) D$ for any fixed $\epsilon>0$,
the probability of a $t$-set being uncovered is
$P(\tau) = O \bigl( (\log v)^{-1/D} \bigr)$.
Thus:

\begin{corollary}
\label{cor:error}
The expected density of a covering produced by the random greedy algorithm is
$1 + O \bigl( (\log v)^{-1/D} \bigr)$, where $D = \textbc{k}{t} - 1$.
\end{corollary}

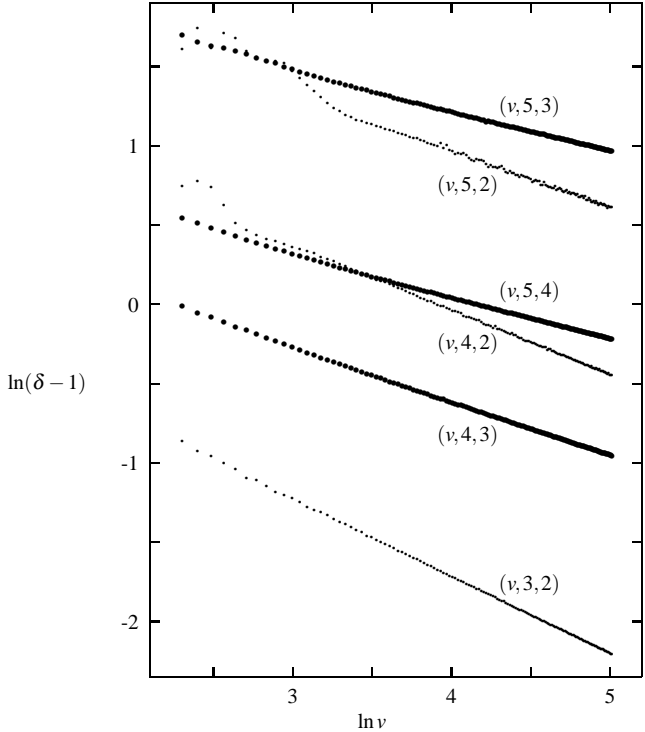
\begin{figure}[t]
\begin{center}
\setlength{\unitlength}{60pt}
\begin{picture}(3.3,4.1)(1.6,-2.4)
\newcommand{\labelsize}{\footnotesize}
\newcommand{\hashmarklength}{.06}
\newcommand{\ttwo}{\circle*{.02}}
\newcommand{\tthree}{\circle*{.03}}
\newcommand{\tfour}{\tthree}
\put(4.3,-1.77){\makebox(0,0)[l]{\labelsize $(v,3,2)$}}
\put(4.3,-0.25){\makebox(0,0)[r]{\labelsize $(v,4,2)$}}
\put(4.3,-0.83){\makebox(0,0)[r]{\labelsize $(v,4,3)$}}
\put(4.3,0.75){\makebox(0,0)[r]{\labelsize $(v,5,2)$}}
\put(4.3,1.25){\makebox(0,0)[l]{\labelsize $(v,5,3)$}}
\put(4.3,0.08){\makebox(0,0)[l]{\labelsize $(v,5,4)$}}
\multiput(2.1,-2.35)(0,4.25){2}{\line(1,0){3.1}}
\multiput(2.1,-2.35)(3.1,0){2}{\line(0,1){4.25}}
\multiput(2.5,-2.35)(.5,0){6}{\line(0,1){\hashmarklength}}
\multiput(2.5,1.9)(.5,0){6}{\line(0,-1){\hashmarklength}}
\multiput(2.1,-2)(0,.5){8}{\line(1,0){\hashmarklength}}
\multiput(5.2,-2)(0,.5){8}{\line(-1,0){\hashmarklength}}
\put(2.03,-2){\makebox(0,0)[r]{\labelsize -2}}
\put(2.03,-1){\makebox(0,0)[r]{\labelsize -1}}
\put(2.03,0){\makebox(0,0)[r]{\labelsize 0}}
\put(2.03,1){\makebox(0,0)[r]{\labelsize 1}}
\put(3,-2.42){\makebox(0,0)[t]{\labelsize 3}}
\put(4,-2.42){\makebox(0,0)[t]{\labelsize 4}}
\put(5,-2.42){\makebox(0,0)[t]{\labelsize 5}}
\put(1.7,-.5){\makebox(0,0)[r]{\labelsize $\ln (\delta-1)$}}
\put(3.5,-2.58){\makebox(0,0)[t]{\labelsize $\ln v$}}
%(v,3,2), 1000x for v in 10..150
\put(2.302585,-0.862118){\ttwo} \put(2.397895,-0.921303){\ttwo}
\put(2.484907,-0.957586){\ttwo} \put(2.564949,-1.002393){\ttwo}
\put(2.639057,-1.038210){\ttwo} \put(2.708050,-1.092857){\ttwo}
\put(2.772589,-1.108814){\ttwo} \put(2.833213,-1.143302){\ttwo}
\put(2.890372,-1.180844){\ttwo} \put(2.944439,-1.202220){\ttwo}
\put(2.995732,-1.221244){\ttwo} \put(3.044522,-1.249468){\ttwo}
\put(3.091042,-1.275102){\ttwo} \put(3.135494,-1.296215){\ttwo}
\put(3.178054,-1.308528){\ttwo} \put(3.218876,-1.329990){\ttwo}
\put(3.258097,-1.353771){\ttwo} \put(3.295837,-1.369547){\ttwo}
\put(3.332205,-1.386294){\ttwo} \put(3.367296,-1.403847){\ttwo}
\put(3.401197,-1.422358){\ttwo} \put(3.433987,-1.437004){\ttwo}
\put(3.465736,-1.460097){\ttwo} \put(3.496508,-1.470368){\ttwo}
\put(3.526361,-1.484951){\ttwo} \put(3.555348,-1.499130){\ttwo}
\put(3.583519,-1.512657){\ttwo} \put(3.610918,-1.526821){\ttwo}
\put(3.637586,-1.537772){\ttwo} \put(3.663562,-1.556223){\ttwo}
\put(3.688879,-1.561527){\ttwo} \put(3.713572,-1.575449){\ttwo}
\put(3.737670,-1.587213){\ttwo} \put(3.761200,-1.598929){\ttwo}
\put(3.784190,-1.607664){\ttwo} \put(3.806662,-1.619825){\ttwo}
\put(3.828641,-1.630884){\ttwo} \put(3.850148,-1.641357){\ttwo}
\put(3.871201,-1.648599){\ttwo} \put(3.891820,-1.662344){\ttwo}
\put(3.912023,-1.673582){\ttwo} \put(3.931826,-1.683300){\ttwo}
\put(3.951244,-1.692180){\ttwo} \put(3.970292,-1.701550){\ttwo}
\put(3.988984,-1.712588){\ttwo} \put(4.007333,-1.717451){\ttwo}
\put(4.025352,-1.726964){\ttwo} \put(4.043051,-1.735393){\ttwo}
\put(4.060443,-1.745181){\ttwo} \put(4.077537,-1.751620){\ttwo}
\put(4.094345,-1.761740){\ttwo} \put(4.110874,-1.769741){\ttwo}
\put(4.127134,-1.776521){\ttwo} \put(4.143135,-1.783078){\ttwo}
\put(4.158883,-1.792939){\ttwo} \put(4.174387,-1.797756){\ttwo}
\put(4.189655,-1.807309){\ttwo} \put(4.204693,-1.813597){\ttwo}
\put(4.219508,-1.820072){\ttwo} \put(4.234107,-1.828036){\ttwo}
\put(4.248495,-1.835209){\ttwo} \put(4.262680,-1.841773){\ttwo}
\put(4.276666,-1.847123){\ttwo} \put(4.290459,-1.856554){\ttwo}
\put(4.304065,-1.860390){\ttwo} \put(4.317488,-1.870227){\ttwo}
\put(4.330733,-1.874295){\ttwo} \put(4.343805,-1.880615){\ttwo}
\put(4.356709,-1.889384){\ttwo} \put(4.369448,-1.894734){\ttwo}
\put(4.382027,-1.901823){\ttwo} \put(4.394449,-1.906996){\ttwo}
\put(4.406719,-1.913563){\ttwo} \put(4.418841,-1.918765){\ttwo}
\put(4.430817,-1.925080){\ttwo} \put(4.442651,-1.928660){\ttwo}
\put(4.454347,-1.933979){\ttwo} \put(4.465908,-1.938911){\ttwo}
\put(4.477337,-1.945470){\ttwo} \put(4.488636,-1.953552){\ttwo}
\put(4.499810,-1.959415){\ttwo} \put(4.510860,-1.963062){\ttwo}
\put(4.521789,-1.969405){\ttwo} \put(4.532599,-1.971880){\ttwo}
\put(4.543295,-1.978455){\ttwo} \put(4.553877,-1.984069){\ttwo}
\put(4.564348,-1.987861){\ttwo} \put(4.574711,-1.993302){\ttwo}
\put(4.584967,-1.997027){\ttwo} \put(4.595120,-2.001926){\ttwo}
\put(4.605170,-2.008822){\ttwo} \put(4.615121,-2.013344){\ttwo}
\put(4.624973,-2.019550){\ttwo} \put(4.634729,-2.022369){\ttwo}
\put(4.644391,-2.029476){\ttwo} \put(4.653960,-2.032784){\ttwo}
\put(4.663439,-2.035364){\ttwo} \put(4.672829,-2.042179){\ttwo}
\put(4.682131,-2.044046){\ttwo} \put(4.691348,-2.049295){\ttwo}
\put(4.700480,-2.055178){\ttwo} \put(4.709530,-2.059633){\ttwo}
\put(4.718499,-2.062064){\ttwo} \put(4.727388,-2.068045){\ttwo}
\put(4.736198,-2.070860){\ttwo} \put(4.744932,-2.076440){\ttwo}
\put(4.753590,-2.077709){\ttwo} \put(4.762174,-2.085629){\ttwo}
\put(4.770685,-2.088259){\ttwo} \put(4.779123,-2.092624){\ttwo}
\put(4.787492,-2.095550){\ttwo} \put(4.795791,-2.100846){\ttwo}
\put(4.804021,-2.104181){\ttwo} \put(4.812184,-2.108083){\ttwo}
\put(4.820282,-2.111341){\ttwo} \put(4.828314,-2.116893){\ttwo}
\put(4.836282,-2.120346){\ttwo} \put(4.844187,-2.124346){\ttwo}
\put(4.852030,-2.128415){\ttwo} \put(4.859812,-2.132634){\ttwo}
\put(4.867534,-2.135889){\ttwo} \put(4.875197,-2.139739){\ttwo}
\put(4.882802,-2.142780){\ttwo} \put(4.890349,-2.146651){\ttwo}
\put(4.897840,-2.150469){\ttwo} \put(4.905275,-2.154217){\ttwo}
\put(4.912655,-2.157252){\ttwo} \put(4.919981,-2.160040){\ttwo}
\put(4.927254,-2.165152){\ttwo} \put(4.934474,-2.169331){\ttwo}
\put(4.941642,-2.170128){\ttwo} \put(4.948760,-2.175736){\ttwo}
\put(4.955827,-2.177500){\ttwo} \put(4.962845,-2.180898){\ttwo}
\put(4.969813,-2.184916){\ttwo} \put(4.976734,-2.188708){\ttwo}
\put(4.983607,-2.191893){\ttwo} \put(4.990433,-2.194316){\ttwo}
\put(4.997212,-2.198265){\ttwo} \put(5.003946,-2.201951){\ttwo}
\put(5.010635,-2.204565){\ttwo}
%(v,4,2), 1000x for v in 10..50, 100x for v in 51..150
\put(2.302585,0.747635){\ttwo} \put(2.397895,0.776820){\ttwo}
\put(2.484907,0.739554){\ttwo} \put(2.564949,0.625115){\ttwo}
\put(2.639057,0.511243){\ttwo} \put(2.708050,0.471253){\ttwo}
\put(2.772589,0.438900){\ttwo} \put(2.833213,0.412012){\ttwo}
\put(2.890372,0.399697){\ttwo} \put(2.944439,0.380957){\ttwo}
\put(2.995732,0.362521){\ttwo} \put(3.044522,0.350556){\ttwo}
\put(3.091042,0.336008){\ttwo} \put(3.135494,0.324201){\ttwo}
\put(3.178054,0.302437){\ttwo} \put(3.218876,0.286606){\ttwo}
\put(3.258097,0.272701){\ttwo} \put(3.295837,0.256958){\ttwo}
\put(3.332205,0.240890){\ttwo} \put(3.367296,0.224738){\ttwo}
\put(3.401197,0.212215){\ttwo} \put(3.433987,0.194346){\ttwo}
\put(3.465736,0.184671){\ttwo} \put(3.496508,0.174908){\ttwo}
\put(3.526361,0.156638){\ttwo} \put(3.555348,0.147637){\ttwo}
\put(3.583519,0.139348){\ttwo} \put(3.610918,0.127505){\ttwo}
\put(3.637586,0.116537){\ttwo} \put(3.663562,0.104105){\ttwo}
\put(3.688879,0.096359){\ttwo} \put(3.713572,0.084038){\ttwo}
\put(3.737670,0.073212){\ttwo} \put(3.761200,0.063489){\ttwo}
\put(3.784190,0.056003){\ttwo} \put(3.806662,0.044510){\ttwo}
\put(3.828641,0.035031){\ttwo} \put(3.850148,0.022212){\ttwo}
\put(3.871201,0.020637){\ttwo} \put(3.891820,0.007599){\ttwo}
\put(3.912023,0.000424){\ttwo} \put(3.931826,-0.010169){\ttwo}
\put(3.951244,-0.011377){\ttwo} \put(3.970292,-0.024834){\ttwo}
\put(3.988984,-0.028494){\ttwo} \put(4.007333,-0.037880){\ttwo}
\put(4.025352,-0.041499){\ttwo} \put(4.043051,-0.049119){\ttwo}
\put(4.060443,-0.052153){\ttwo} \put(4.077537,-0.074351){\ttwo}
\put(4.094345,-0.076588){\ttwo} \put(4.110874,-0.080535){\ttwo}
\put(4.127134,-0.085868){\ttwo} \put(4.143135,-0.088269){\ttwo}
\put(4.158883,-0.095325){\ttwo} \put(4.174387,-0.099980){\ttwo}
\put(4.189655,-0.108473){\ttwo} \put(4.204693,-0.123157){\ttwo}
\put(4.219508,-0.130001){\ttwo} \put(4.234107,-0.131123){\ttwo}
\put(4.248495,-0.137479){\ttwo} \put(4.262680,-0.144665){\ttwo}
\put(4.276666,-0.135680){\ttwo} \put(4.290459,-0.157696){\ttwo}
\put(4.304065,-0.156158){\ttwo} \put(4.317488,-0.161947){\ttwo}
\put(4.330733,-0.168979){\ttwo} \put(4.343805,-0.174248){\ttwo}
\put(4.356709,-0.177498){\ttwo} \put(4.369448,-0.182477){\ttwo}
\put(4.382027,-0.180058){\ttwo} \put(4.394449,-0.193383){\ttwo}
\put(4.406719,-0.192235){\ttwo} \put(4.418841,-0.205124){\ttwo}
\put(4.430817,-0.212444){\ttwo} \put(4.442651,-0.211136){\ttwo}
\put(4.454347,-0.222118){\ttwo} \put(4.465908,-0.230690){\ttwo}
\put(4.477337,-0.229531){\ttwo} \put(4.488636,-0.226309){\ttwo}
\put(4.499810,-0.238428){\ttwo} \put(4.510860,-0.241535){\ttwo}
\put(4.521789,-0.244878){\ttwo} \put(4.532599,-0.252046){\ttwo}
\put(4.543295,-0.256377){\ttwo} \put(4.553877,-0.258446){\ttwo}
\put(4.564348,-0.263161){\ttwo} \put(4.574711,-0.268790){\ttwo}
\put(4.584967,-0.270398){\ttwo} \put(4.595120,-0.270231){\ttwo}
\put(4.605170,-0.285825){\ttwo} \put(4.615121,-0.280998){\ttwo}
\put(4.624973,-0.283227){\ttwo} \put(4.634729,-0.294444){\ttwo}
\put(4.644391,-0.291673){\ttwo} \put(4.653960,-0.298732){\ttwo}
\put(4.663439,-0.304974){\ttwo} \put(4.672829,-0.309562){\ttwo}
\put(4.682131,-0.312423){\ttwo} \put(4.691348,-0.319582){\ttwo}
\put(4.700480,-0.321808){\ttwo} \put(4.709530,-0.320347){\ttwo}
\put(4.718499,-0.326321){\ttwo} \put(4.727388,-0.330187){\ttwo}
\put(4.736198,-0.333589){\ttwo} \put(4.744932,-0.336994){\ttwo}
\put(4.753590,-0.339896){\ttwo} \put(4.762174,-0.341880){\ttwo}
\put(4.770685,-0.347405){\ttwo} \put(4.779123,-0.352634){\ttwo}
\put(4.787492,-0.356855){\ttwo} \put(4.795791,-0.355200){\ttwo}
\put(4.804021,-0.364232){\ttwo} \put(4.812184,-0.365337){\ttwo}
\put(4.820282,-0.367524){\ttwo} \put(4.828314,-0.370465){\ttwo}
\put(4.836282,-0.373634){\ttwo} \put(4.844187,-0.380745){\ttwo}
\put(4.852030,-0.380508){\ttwo} \put(4.859812,-0.382014){\ttwo}
\put(4.867534,-0.385136){\ttwo} \put(4.875197,-0.391482){\ttwo}
\put(4.882802,-0.390604){\ttwo} \put(4.890349,-0.399162){\ttwo}
\put(4.897840,-0.401662){\ttwo} \put(4.905275,-0.403477){\ttwo}
\put(4.912655,-0.401942){\ttwo} \put(4.919981,-0.407281){\ttwo}
\put(4.927254,-0.410987){\ttwo} \put(4.934474,-0.410372){\ttwo}
\put(4.941642,-0.418525){\ttwo} \put(4.948760,-0.420779){\ttwo}
\put(4.955827,-0.420349){\ttwo} \put(4.962845,-0.428989){\ttwo}
\put(4.969813,-0.427918){\ttwo} \put(4.976734,-0.430916){\ttwo}
\put(4.983607,-0.435569){\ttwo} \put(4.990433,-0.432526){\ttwo}
\put(4.997212,-0.441618){\ttwo} \put(5.003946,-0.445542){\ttwo}
\put(5.010635,-0.445742){\ttwo}
%(v,4,3), 1000x for v in 10..50, 100x for v in 51..150
\put(2.302585,-0.008234){\tthree} \put(2.397895,-0.051868){\tthree}
\put(2.484907,-0.076214){\tthree} \put(2.564949,-0.108723){\tthree}
\put(2.639057,-0.138921){\tthree} \put(2.708050,-0.162369){\tthree}
\put(2.772589,-0.189001){\tthree} \put(2.833213,-0.209741){\tthree}
\put(2.890372,-0.229228){\tthree} \put(2.944439,-0.251658){\tthree}
\put(2.995732,-0.270401){\tthree} \put(3.044522,-0.285675){\tthree}
\put(3.091042,-0.304456){\tthree} \put(3.135494,-0.320434){\tthree}
\put(3.178054,-0.335478){\tthree} \put(3.218876,-0.349062){\tthree}
\put(3.258097,-0.362054){\tthree} \put(3.295837,-0.377461){\tthree}
\put(3.332205,-0.388619){\tthree} \put(3.367296,-0.400718){\tthree}
\put(3.401197,-0.411515){\tthree} \put(3.433987,-0.424450){\tthree}
\put(3.465736,-0.435116){\tthree} \put(3.496508,-0.446004){\tthree}
\put(3.526361,-0.456534){\tthree} \put(3.555348,-0.466707){\tthree}
\put(3.583519,-0.476966){\tthree} \put(3.610918,-0.485517){\tthree}
\put(3.637586,-0.494308){\tthree} \put(3.663562,-0.502898){\tthree}
\put(3.688879,-0.511741){\tthree} \put(3.713572,-0.520436){\tthree}
\put(3.737670,-0.527816){\tthree} \put(3.761200,-0.536382){\tthree}
\put(3.784190,-0.543784){\tthree} \put(3.806662,-0.551601){\tthree}
\put(3.828641,-0.559474){\tthree} \put(3.850148,-0.566300){\tthree}
\put(3.871201,-0.574231){\tthree} \put(3.891820,-0.581376){\tthree}
\put(3.912023,-0.588139){\tthree} \put(3.931826,-0.594207){\tthree}
\put(3.951244,-0.599266){\tthree} \put(3.970292,-0.607233){\tthree}
\put(3.988984,-0.614244){\tthree} \put(4.007333,-0.620069){\tthree}
\put(4.025352,-0.625917){\tthree} \put(4.043051,-0.632086){\tthree}
\put(4.060443,-0.636175){\tthree} \put(4.077537,-0.643169){\tthree}
\put(4.094345,-0.648195){\tthree} \put(4.110874,-0.653328){\tthree}
\put(4.127134,-0.661040){\tthree} \put(4.143135,-0.667035){\tthree}
\put(4.158883,-0.671308){\tthree} \put(4.174387,-0.676101){\tthree}
\put(4.189655,-0.681170){\tthree} \put(4.204693,-0.685328){\tthree}
\put(4.219508,-0.689556){\tthree} \put(4.234107,-0.694171){\tthree}
\put(4.248495,-0.701930){\tthree} \put(4.262680,-0.704859){\tthree}
\put(4.276666,-0.709330){\tthree} \put(4.290459,-0.714693){\tthree}
\put(4.304065,-0.719356){\tthree} \put(4.317488,-0.724063){\tthree}
\put(4.330733,-0.728759){\tthree} \put(4.343805,-0.731795){\tthree}
\put(4.356709,-0.736335){\tthree} \put(4.369448,-0.741678){\tthree}
\put(4.382027,-0.746542){\tthree} \put(4.394449,-0.749231){\tthree}
\put(4.406719,-0.753870){\tthree} \put(4.418841,-0.757845){\tthree}
\put(4.430817,-0.761724){\tthree} \put(4.442651,-0.764443){\tthree}
\put(4.454347,-0.769519){\tthree} \put(4.465908,-0.772450){\tthree}
\put(4.477337,-0.778548){\tthree} \put(4.488636,-0.780229){\tthree}
\put(4.499810,-0.784403){\tthree} \put(4.510860,-0.788189){\tthree}
\put(4.521789,-0.791194){\tthree} \put(4.532599,-0.794446){\tthree}
\put(4.543295,-0.798892){\tthree} \put(4.553877,-0.801769){\tthree}
\put(4.564348,-0.804905){\tthree} \put(4.574711,-0.809006){\tthree}
\put(4.584967,-0.812599){\tthree} \put(4.595120,-0.815854){\tthree}
\put(4.605170,-0.819575){\tthree} \put(4.615121,-0.821924){\tthree}
\put(4.624973,-0.825822){\tthree} \put(4.634729,-0.829602){\tthree}
\put(4.644391,-0.833158){\tthree} \put(4.653960,-0.835002){\tthree}
\put(4.663439,-0.837589){\tthree} \put(4.672829,-0.841496){\tthree}
\put(4.682131,-0.844135){\tthree} \put(4.691348,-0.847348){\tthree}
\put(4.700480,-0.850649){\tthree} \put(4.709530,-0.853652){\tthree}
\put(4.718499,-0.857300){\tthree} \put(4.727388,-0.859433){\tthree}
\put(4.736198,-0.862218){\tthree} \put(4.744932,-0.865293){\tthree}
\put(4.753590,-0.867702){\tthree} \put(4.762174,-0.870700){\tthree}
\put(4.770685,-0.874576){\tthree} \put(4.779123,-0.877071){\tthree}
\put(4.787492,-0.879849){\tthree} \put(4.795791,-0.883327){\tthree}
\put(4.804021,-0.884873){\tthree} \put(4.812184,-0.887438){\tthree}
\put(4.820282,-0.890564){\tthree} \put(4.828314,-0.893554){\tthree}
\put(4.836282,-0.896185){\tthree} \put(4.844187,-0.898577){\tthree}
\put(4.852030,-0.900746){\tthree} \put(4.859812,-0.903833){\tthree}
\put(4.867534,-0.906411){\tthree} \put(4.875197,-0.908639){\tthree}
\put(4.882802,-0.911459){\tthree} \put(4.890349,-0.913367){\tthree}
\put(4.897840,-0.916263){\tthree} \put(4.905275,-0.918771){\tthree}
\put(4.912655,-0.920676){\tthree} \put(4.919981,-0.923313){\tthree}
\put(4.927254,-0.925889){\tthree} \put(4.934474,-0.928457){\tthree}
\put(4.941642,-0.930639){\tthree} \put(4.948760,-0.933209){\tthree}
\put(4.955827,-0.935625){\tthree} \put(4.962845,-0.937621){\tthree}
\put(4.969813,-0.939467){\tthree} \put(4.976734,-0.942815){\tthree}
\put(4.983607,-0.944399){\tthree} \put(4.990433,-0.946746){\tthree}
\put(4.997212,-0.949301){\tthree} \put(5.003946,-0.951365){\tthree}
\put(5.010635,-0.953685){\tthree}
%(v,5,2), 1000x for v in 10..50, 10x for v in 51..150
\put(2.302585,1.609438){\ttwo} \put(2.397895,1.745239){\ttwo}
\put(2.484907,1.618217){\ttwo} \put(2.564949,1.711717){\ttwo}
\put(2.639057,1.682913){\ttwo} \put(2.708050,1.597481){\ttwo}
\put(2.772589,1.562661){\ttwo} \put(2.833213,1.539615){\ttwo}
\put(2.890372,1.535855){\ttwo} \put(2.944439,1.521814){\ttwo}
\put(2.995732,1.472858){\ttwo} \put(3.044522,1.429969){\ttwo}
\put(3.091042,1.385645){\ttwo} \put(3.135494,1.345781){\ttwo}
\put(3.178054,1.308450){\ttwo} \put(3.218876,1.272902){\ttwo}
\put(3.258097,1.240352){\ttwo} \put(3.295837,1.219165){\ttwo}
\put(3.332205,1.194644){\ttwo} \put(3.367296,1.182867){\ttwo}
\put(3.401197,1.165820){\ttwo} \put(3.433987,1.157611){\ttwo}
\put(3.465736,1.144229){\ttwo} \put(3.496508,1.135974){\ttwo}
\put(3.526361,1.126664){\ttwo} \put(3.555348,1.118761){\ttwo}
\put(3.583519,1.109128){\ttwo} \put(3.610918,1.100457){\ttwo}
\put(3.637586,1.091747){\ttwo} \put(3.663562,1.083766){\ttwo}
\put(3.688879,1.074241){\ttwo} \put(3.713572,1.067839){\ttwo}
\put(3.737670,1.059567){\ttwo} \put(3.761200,1.049764){\ttwo}
\put(3.784190,1.047430){\ttwo} \put(3.806662,1.034704){\ttwo}
\put(3.828641,1.029371){\ttwo} \put(3.850148,1.019365){\ttwo}
\put(3.871201,1.012741){\ttwo} \put(3.891820,1.007073){\ttwo}
\put(3.912023,1.000058){\ttwo} \put(3.931826,1.017785){\ttwo}
\put(3.951244,0.986188){\ttwo} \put(3.970292,1.004318){\ttwo}
\put(3.988984,0.976101){\ttwo} \put(4.007333,0.969401){\ttwo}
\put(4.025352,0.951257){\ttwo} \put(4.043051,0.955367){\ttwo}
\put(4.060443,0.948085){\ttwo} \put(4.077537,0.956949){\ttwo}
\put(4.094345,0.943048){\ttwo} \put(4.110874,0.923478){\ttwo}
\put(4.127134,0.910031){\ttwo} \put(4.143135,0.904237){\ttwo}
\put(4.158883,0.916291){\ttwo} \put(4.174387,0.907987){\ttwo}
\put(4.189655,0.909836){\ttwo} \put(4.204693,0.907843){\ttwo}
\put(4.219508,0.889418){\ttwo} \put(4.234107,0.869698){\ttwo}
\put(4.248495,0.883203){\ttwo} \put(4.262680,0.865357){\ttwo}
\put(4.276666,0.864256){\ttwo} \put(4.290459,0.862832){\ttwo}
\put(4.304065,0.865329){\ttwo} \put(4.317488,0.872762){\ttwo}
\put(4.330733,0.855535){\ttwo} \put(4.343805,0.834572){\ttwo}
\put(4.356709,0.829733){\ttwo} \put(4.369448,0.845070){\ttwo}
\put(4.382027,0.827252){\ttwo} \put(4.394449,0.832909){\ttwo}
\put(4.406719,0.819027){\ttwo} \put(4.418841,0.817082){\ttwo}
\put(4.430817,0.802287){\ttwo} \put(4.442651,0.815960){\ttwo}
\put(4.454347,0.812297){\ttwo} \put(4.465908,0.802070){\ttwo}
\put(4.477337,0.804641){\ttwo} \put(4.488636,0.798597){\ttwo}
\put(4.499810,0.793552){\ttwo} \put(4.510860,0.774373){\ttwo}
\put(4.521789,0.778395){\ttwo} \put(4.532599,0.783921){\ttwo}
\put(4.543295,0.781445){\ttwo} \put(4.553877,0.778842){\ttwo}
\put(4.564348,0.754081){\ttwo} \put(4.574711,0.754378){\ttwo}
\put(4.584967,0.772105){\ttwo} \put(4.595120,0.753444){\ttwo}
\put(4.605170,0.747119){\ttwo} \put(4.615121,0.742503){\ttwo}
\put(4.624973,0.742021){\ttwo} \put(4.634729,0.737276){\ttwo}
\put(4.644391,0.735730){\ttwo} \put(4.653960,0.748977){\ttwo}
\put(4.663439,0.725241){\ttwo} \put(4.672829,0.728828){\ttwo}
\put(4.682131,0.732098){\ttwo} \put(4.691348,0.730661){\ttwo}
\put(4.700480,0.715253){\ttwo} \put(4.709530,0.711648){\ttwo}
\put(4.718499,0.717227){\ttwo} \put(4.727388,0.723804){\ttwo}
\put(4.736198,0.701958){\ttwo} \put(4.744932,0.707686){\ttwo}
\put(4.753590,0.701062){\ttwo} \put(4.762174,0.705013){\ttwo}
\put(4.770685,0.701586){\ttwo} \put(4.779123,0.694286){\ttwo}
\put(4.787492,0.689569){\ttwo} \put(4.795791,0.694180){\ttwo}
\put(4.804021,0.677650){\ttwo} \put(4.812184,0.677227){\ttwo}
\put(4.820282,0.683264){\ttwo} \put(4.828314,0.686025){\ttwo}
\put(4.836282,0.673527){\ttwo} \put(4.844187,0.682404){\ttwo}
\put(4.852030,0.658091){\ttwo} \put(4.859812,0.657829){\ttwo}
\put(4.867534,0.671812){\ttwo} \put(4.875197,0.662985){\ttwo}
\put(4.882802,0.644990){\ttwo} \put(4.890349,0.648237){\ttwo}
\put(4.897840,0.649157){\ttwo} \put(4.905275,0.638969){\ttwo}
\put(4.912655,0.652632){\ttwo} \put(4.919981,0.640474){\ttwo}
\put(4.927254,0.628418){\ttwo} \put(4.934474,0.635032){\ttwo}
\put(4.941642,0.639850){\ttwo} \put(4.948760,0.633823){\ttwo}
\put(4.955827,0.629747){\ttwo} \put(4.962845,0.636516){\ttwo}
\put(4.969813,0.627209){\ttwo} \put(4.976734,0.622151){\ttwo}
\put(4.983607,0.612053){\ttwo} \put(4.990433,0.620694){\ttwo}
\put(4.997212,0.616512){\ttwo} \put(5.003946,0.614543){\ttwo}
\put(5.010635,0.618192){\ttwo}
%(v,5,3), 1000x for v in 10..50, 10x for v in 51..150
\put(2.302585,1.697449){\tthree} \put(2.397895,1.653855){\tthree}
\put(2.484907,1.632631){\tthree} \put(2.564949,1.615969){\tthree}
\put(2.639057,1.601982){\tthree} \put(2.708050,1.578145){\tthree}
\put(2.772589,1.557550){\tthree} \put(2.833213,1.536231){\tthree}
\put(2.890372,1.520341){\tthree} \put(2.944439,1.500024){\tthree}
\put(2.995732,1.483230){\tthree} \put(3.044522,1.468131){\tthree}
\put(3.091042,1.454655){\tthree} \put(3.135494,1.440350){\tthree}
\put(3.178054,1.427866){\tthree} \put(3.218876,1.416813){\tthree}
\put(3.258097,1.406795){\tthree} \put(3.295837,1.396435){\tthree}
\put(3.332205,1.386157){\tthree} \put(3.367296,1.376269){\tthree}
\put(3.401197,1.367546){\tthree} \put(3.433987,1.359572){\tthree}
\put(3.465736,1.350014){\tthree} \put(3.496508,1.342322){\tthree}
\put(3.526361,1.334000){\tthree} \put(3.555348,1.326813){\tthree}
\put(3.583519,1.320666){\tthree} \put(3.610918,1.311988){\tthree}
\put(3.637586,1.305337){\tthree} \put(3.663562,1.299356){\tthree}
\put(3.688879,1.292416){\tthree} \put(3.713572,1.286376){\tthree}
\put(3.737670,1.280611){\tthree} \put(3.761200,1.273420){\tthree}
\put(3.784190,1.268099){\tthree} \put(3.806662,1.262656){\tthree}
\put(3.828641,1.256359){\tthree} \put(3.850148,1.251153){\tthree}
\put(3.871201,1.246008){\tthree} \put(3.891820,1.241107){\tthree}
\put(3.912023,1.235570){\tthree} \put(3.931826,1.234117){\tthree}
\put(3.951244,1.226779){\tthree} \put(3.970292,1.222665){\tthree}
\put(3.988984,1.218060){\tthree} \put(4.007333,1.212411){\tthree}
\put(4.025352,1.207894){\tthree} \put(4.043051,1.202861){\tthree}
\put(4.060443,1.196670){\tthree} \put(4.077537,1.194680){\tthree}
\put(4.094345,1.190063){\tthree} \put(4.110874,1.183844){\tthree}
\put(4.127134,1.182220){\tthree} \put(4.143135,1.174596){\tthree}
\put(4.158883,1.175771){\tthree} \put(4.174387,1.168630){\tthree}
\put(4.189655,1.165838){\tthree} \put(4.204693,1.163568){\tthree}
\put(4.219508,1.153960){\tthree} \put(4.234107,1.155500){\tthree}
\put(4.248495,1.153439){\tthree} \put(4.262680,1.148988){\tthree}
\put(4.276666,1.146348){\tthree} \put(4.290459,1.142216){\tthree}
\put(4.304065,1.140663){\tthree} \put(4.317488,1.134467){\tthree}
\put(4.330733,1.133900){\tthree} \put(4.343805,1.128098){\tthree}
\put(4.356709,1.127657){\tthree} \put(4.369448,1.122021){\tthree}
\put(4.382027,1.120510){\tthree} \put(4.394449,1.119042){\tthree}
\put(4.406719,1.114342){\tthree} \put(4.418841,1.112651){\tthree}
\put(4.430817,1.107878){\tthree} \put(4.442651,1.107252){\tthree}
\put(4.454347,1.102902){\tthree} \put(4.465908,1.099951){\tthree}
\put(4.477337,1.096761){\tthree} \put(4.488636,1.093722){\tthree}
\put(4.499810,1.090356){\tthree} \put(4.510860,1.088935){\tthree}
\put(4.521789,1.087297){\tthree} \put(4.532599,1.081066){\tthree}
\put(4.543295,1.080606){\tthree} \put(4.553877,1.076737){\tthree}
\put(4.564348,1.075639){\tthree} \put(4.574711,1.075004){\tthree}
\put(4.584967,1.073169){\tthree} \put(4.595120,1.070605){\tthree}
\put(4.605170,1.065625){\tthree} \put(4.615121,1.064276){\tthree}
\put(4.624973,1.061333){\tthree} \put(4.634729,1.059762){\tthree}
\put(4.644391,1.056498){\tthree} \put(4.653960,1.053343){\tthree}
\put(4.663439,1.050598){\tthree} \put(4.672829,1.048137){\tthree}
\put(4.682131,1.048287){\tthree} \put(4.691348,1.044591){\tthree}
\put(4.700480,1.043742){\tthree} \put(4.709530,1.041888){\tthree}
\put(4.718499,1.038836){\tthree} \put(4.727388,1.036701){\tthree}
\put(4.736198,1.034752){\tthree} \put(4.744932,1.031663){\tthree}
\put(4.753590,1.030365){\tthree} \put(4.762174,1.028663){\tthree}
\put(4.770685,1.025256){\tthree} \put(4.779123,1.024321){\tthree}
\put(4.787492,1.021697){\tthree} \put(4.795791,1.019695){\tthree}
\put(4.804021,1.017236){\tthree} \put(4.812184,1.015481){\tthree}
\put(4.820282,1.013898){\tthree} \put(4.828314,1.010990){\tthree}
\put(4.836282,1.010441){\tthree} \put(4.844187,1.008212){\tthree}
\put(4.852030,1.006706){\tthree} \put(4.859812,1.004705){\tthree}
\put(4.867534,1.003449){\tthree} \put(4.875197,1.000784){\tthree}
\put(4.882802,0.999893){\tthree} \put(4.890349,0.997029){\tthree}
\put(4.897840,0.995543){\tthree} \put(4.905275,0.993281){\tthree}
\put(4.912655,0.992642){\tthree} \put(4.919981,0.989548){\tthree}
\put(4.927254,0.987903){\tthree} \put(4.934474,0.986956){\tthree}
\put(4.941642,0.985048){\tthree} \put(4.948760,0.983724){\tthree}
\put(4.955827,0.981041){\tthree} \put(4.962845,0.979790){\tthree}
\put(4.969813,0.977574){\tthree} \put(4.976734,0.977535){\tthree}
\put(4.983607,0.974410){\tthree} \put(4.990433,0.973288){\tthree}
\put(4.997212,0.970731){\tthree} \put(5.003946,0.969457){\tthree}
\put(5.010635,0.968608){\tthree}
%(v,5,4), 1000x for v in 10..50, 10x for v in 51..150
\put(2.302585,0.545862){\tfour} \put(2.397895,0.512860){\tfour}
\put(2.484907,0.480904){\tfour} \put(2.564949,0.454499){\tfour}
\put(2.639057,0.431510){\tfour} \put(2.708050,0.407544){\tfour}
\put(2.772589,0.387264){\tfour} \put(2.833213,0.368912){\tfour}
\put(2.890372,0.352207){\tfour} \put(2.944439,0.334746){\tfour}
\put(2.995732,0.319214){\tfour} \put(3.044522,0.305588){\tfour}
\put(3.091042,0.291161){\tfour} \put(3.135494,0.277957){\tfour}
\put(3.178054,0.266724){\tfour} \put(3.218876,0.255191){\tfour}
\put(3.258097,0.243118){\tfour} \put(3.295837,0.232395){\tfour}
\put(3.332205,0.222627){\tfour} \put(3.367296,0.212664){\tfour}
\put(3.401197,0.203696){\tfour} \put(3.433987,0.194206){\tfour}
\put(3.465736,0.185453){\tfour} \put(3.496508,0.176998){\tfour}
\put(3.526361,0.169053){\tfour} \put(3.555348,0.161128){\tfour}
\put(3.583519,0.153301){\tfour} \put(3.610918,0.145895){\tfour}
\put(3.637586,0.138693){\tfour} \put(3.663562,0.131806){\tfour}
\put(3.688879,0.125026){\tfour} \put(3.713572,0.118379){\tfour}
\put(3.737670,0.111688){\tfour} \put(3.761200,0.105497){\tfour}
\put(3.784190,0.099353){\tfour} \put(3.806662,0.093068){\tfour}
\put(3.828641,0.087559){\tfour} \put(3.850148,0.081708){\tfour}
\put(3.871201,0.076338){\tfour} \put(3.891820,0.070853){\tfour}
\put(3.912023,0.065469){\tfour} \put(3.931826,0.059627){\tfour}
\put(3.951244,0.055972){\tfour} \put(3.970292,0.048958){\tfour}
\put(3.988984,0.045129){\tfour} \put(4.007333,0.040164){\tfour}
\put(4.025352,0.036627){\tfour} \put(4.043051,0.030688){\tfour}
\put(4.060443,0.026093){\tfour} \put(4.077537,0.020813){\tfour}
\put(4.094345,0.017327){\tfour} \put(4.110874,0.012857){\tfour}
\put(4.127134,0.009017){\tfour} \put(4.143135,0.004812){\tfour}
\put(4.158883,0.001400){\tfour} \put(4.174387,-0.003041){\tfour}
\put(4.189655,-0.007073){\tfour} \put(4.204693,-0.011004){\tfour}
\put(4.219508,-0.014542){\tfour} \put(4.234107,-0.018612){\tfour}
\put(4.248495,-0.022798){\tfour} \put(4.262680,-0.026649){\tfour}
\put(4.276666,-0.029975){\tfour} \put(4.290459,-0.033224){\tfour}
\put(4.304065,-0.036927){\tfour} \put(4.317488,-0.041043){\tfour}
\put(4.330733,-0.043884){\tfour} \put(4.343805,-0.046931){\tfour}
\put(4.356709,-0.050583){\tfour} \put(4.369448,-0.053706){\tfour}
\put(4.382027,-0.056479){\tfour} \put(4.394449,-0.060192){\tfour}
\put(4.406719,-0.063128){\tfour} \put(4.418841,-0.065839){\tfour}
\put(4.430817,-0.069360){\tfour} \put(4.442651,-0.072355){\tfour}
\put(4.454347,-0.075132){\tfour} \put(4.465908,-0.078269){\tfour}
\put(4.477337,-0.081085){\tfour} \put(4.488636,-0.083821){\tfour}
\put(4.499810,-0.087027){\tfour} \put(4.510860,-0.090449){\tfour}
\put(4.521789,-0.092944){\tfour} \put(4.532599,-0.095569){\tfour}
\put(4.543295,-0.097822){\tfour} \put(4.553877,-0.101196){\tfour}
\put(4.564348,-0.103396){\tfour} \put(4.574711,-0.107046){\tfour}
\put(4.584967,-0.108833){\tfour} \put(4.595120,-0.111005){\tfour}
\put(4.605170,-0.114213){\tfour} \put(4.615121,-0.116829){\tfour}
\put(4.624973,-0.119014){\tfour} \put(4.634729,-0.121266){\tfour}
\put(4.644391,-0.124035){\tfour} \put(4.653960,-0.126621){\tfour}
\put(4.663439,-0.128879){\tfour} \put(4.672829,-0.131380){\tfour}
\put(4.682131,-0.133923){\tfour} \put(4.691348,-0.135742){\tfour}
\put(4.700480,-0.138655){\tfour} \put(4.709530,-0.141091){\tfour}
\put(4.718499,-0.143144){\tfour} \put(4.727388,-0.145659){\tfour}
\put(4.736198,-0.147592){\tfour} \put(4.744932,-0.150137){\tfour}
\put(4.753590,-0.152206){\tfour} \put(4.762174,-0.154278){\tfour}
\put(4.770685,-0.156216){\tfour} \put(4.779123,-0.158451){\tfour}
\put(4.787492,-0.160560){\tfour} \put(4.795791,-0.162979){\tfour}
\put(4.804021,-0.164981){\tfour} \put(4.812184,-0.167082){\tfour}
\put(4.820282,-0.169381){\tfour} \put(4.828314,-0.171689){\tfour}
\put(4.836282,-0.173316){\tfour} \put(4.844187,-0.175231){\tfour}
\put(4.852030,-0.177167){\tfour} \put(4.859812,-0.178737){\tfour}
\put(4.867534,-0.181022){\tfour} \put(4.875197,-0.182929){\tfour}
\put(4.882802,-0.184951){\tfour} \put(4.890349,-0.186926){\tfour}
\put(4.897840,-0.188973){\tfour} \put(4.905275,-0.191014){\tfour}
\put(4.912655,-0.192611){\tfour} \put(4.919981,-0.194675){\tfour}
\put(4.927254,-0.196713){\tfour} \put(4.934474,-0.198571){\tfour}
\put(4.941642,-0.200267){\tfour} \put(4.948760,-0.202163){\tfour}
\put(4.955827,-0.204203){\tfour} \put(4.962845,-0.205921){\tfour}
\put(4.969813,-0.207406){\tfour} \put(4.976734,-0.209177){\tfour}
\put(4.983607,-0.210930){\tfour} \put(4.990433,-0.212468){\tfour}
\put(4.997212,-0.214292){\tfour} \put(5.003946,-0.216229){\tfour}
\put(5.010635,-0.218548){\tfour}
\end{picture}
\end{center}
\caption{Average density $\delta$ of random greedy coverings.}
\label{log-log}
\end{figure}

This bound is pessimistic. Figure~\ref{log-log} gives log-log plots for several
$(k,t)$ pairs, based on 1000 random greedy coverings per $(v,k,t)$ triple for
$v \leq 50$, and $10^{6-k}$ such coverings for $v > 50$. The apparent
asymptotic linearity of the plots suggests that the expected density of a
random greedy covering for $k$~and~$t$ fixed is $1 + \Theta(v^{-\alpha})$, for
some positive $\alpha = \alpha(k,t)$ as $v \to \infty$.

\begin{table}[t]
\vspace*{-2ex}
\begin{center}
\begin{tabular}{|c|c|c|c|}
\hline
$k$ & $t$ &  $\alpha$ & $(k-t)/D$ \\
\hline
3 & 2 & 0.484 & 1/2 \\
4 & 2 & 0.407 & 2/5 \\
4 & 3 & 0.332 & 1/3 \\
5 & 2 & 0.344 & 1/3 \\
5 & 3 & 0.241 & 2/9 \\
5 & 4 & 0.256 & 1/4 \\
\hline
\end{tabular}
\end{center}
\caption{Estimates for $\alpha(k,t)$}\label{alpha-tab}
\end{table}

To estimate~$\alpha$ for each of the curves in Figure~\ref{log-log}, we used
the tails of the curves ($100 \leq v \leq 150$) for a least-squares fit to a
straight line. That gave us rough estimates for the slopes $-\alpha(k,t)$, as
indicated in Table~\ref{alpha-tab}. Those values suggest:

\begin{conjecture}
The expected density of a covering produced by the random greedy algorithm is
$1 + \Theta \bigl( v^{-(k-t)/D} \bigr)$, where $D = \textbc{k}{t} - 1$.
\end{conjecture}

The following argument, though far from a proof, supports the conjecture.

\begin{proof}[Heuristic argument]
Let $\alpha = (k-t)/ \muspace{-2} D$. The conjecture is equivalent to the
statement that there are $\Theta(v^{t-\alpha})$ expected $t$-sets not covered
by a random greedy packing. (The first three steps of Algorithm~\ref{a:greedy}
constitute the random greedy packing algorithm.) So consider the $t$-uniform
hypergraph whose edges are the $t$-sets still uncovered during the packing
algorithm. Assume that this hypergraph looks like a random hypergraph with the
same number of edges, and assume that the packing algorithm has managed to
leave just $c_1 v^{t-\alpha}(1+o(1))$ edges in the hypergraph, for some
positive constant~$c_1$. We show that a positive fraction of these edges---%
that is, $\Theta(v^{t-\alpha})$ in all---% are contained in no $k$-clique, and
hence can never be covered by the packing; this provides the
$\Omega(v^{-\alpha})$ lower bound of the conjecture's error term.

Under the stated assumptions, the probability~$p$ that a given edge exists in
the hypergraph is asymptotic to $c_1 t! \muspace{2} v^{-\alpha} = c_2
v^{-\alpha} \!$, and the probability, for a given edge in the hypergraph and a
given $k$-set containing that edge, that the other $\textbc{k}{t} - 1 = D$
edges on those $k$~vertices also exist is~$p^D \!$. Therefore the expected
number of $k$-cliques that contain the given edge is asymptotic to $p^D
v^{k-t}\!/(k-t)! = c_3 v^{-\alpha D} v^{k-t} = c_3$, a positive constant. But
this number of $k$-cliques is Poisson distributed, so is zero with probability
asymptotic to $e^{-c_3} \!$, also a positive constant, thus a positive fraction
of the edges are contained in no $k$-clique, as claimed. The matching
$O(v^{-\alpha})$ upper bound follows from similar reasoning, and that completes
the argument. It, together with our empirical data, makes the conjecture quite
compelling.
\end{proof}

\section{Induced Coverings and R\"odl's Bound}
\label{s:induced}

While the greedy algorithm produces good coverings, it works in time and space
$\Theta(v^k)$. These can be reduced to time $O(v^t \log v)$ and space $O(v^t)$
using the early abort strategy of Corollary~\ref{cor:abort}, but for larger
values of $v$,~$k$, and~$t$, the induced covering algorithm is more practical,
because it is faster.

\begin{theorem}
\label{t:induced}
For fixed $k$~and~$t$ the expected density
of an induced covering is $1 + o(1)$.
\end{theorem}

\begin{proof}
For Step~\ref{induced-choose} of Algorithm~\ref{a:induced}
choose $\ell = \frac{1}{9} v^{1-1/t} \!$, and choose the prime~$p$ such that
$$
4\ell \drel{\leq} \frac{v-t}{p} \drel{\leq} 8\ell \,.
$$
Such a prime exists by Bertrand's Postulate, which states that there is always
a prime between $n$~and~$2n$. These choices ensure that $p^t>v$, and  that the
affine $(p^t,p^{t-1},t)$ covering by hyperplanes has density $1 + O(v^{-1/t})$.

By Corollary~\ref{cor:error} the precomputed $(\ell',k,t)$ greedy coverings of
Step~\ref{induced-precompute} have expected density $1+O((\log v)^{-1/D})$. So
by running $O(\log\ell)$ trials per precomputed covering, we can ensure, with
probability greater than, for example, $1-1/\ell$, that all precomputed
coverings have density $1+O((\log v)^{-1/D})$.

Now select the $v$-set~$V$ as a random subset of the points in the affine
covering, and consider a fixed $t$-set~$T$ of~$V$. There are, on average, $1 +
O(v^{-1/t})$ hyperplanes containing~$T$; let $P$ be one of them. The size of
the intersection of $V$ and~$P-T$ has a hypergeometric distribution from
$0$~to~$v-t$ with mean
$$M = \frac{(v-t)(p^{t-1}-t)}{p^t-t}.$$
For $p \geq 5$ we have
$$(v-t)/2p < M < (v-t)/p,$$
thus $2\ell < M < 8\ell$ by our choice of~$p$. So the
probability that the size of the intersection is at most~$\ell$ or at
least~$9\ell$ is $O(e^{-c\ell})$ for some $c>0$.

This intersection, together with $T$~itself, is replaced in the induced
covering by an $(\ell',k,t)$ covering. If $\ell < \ell' < 9\ell$, then this
covering has density $1+O((\log v)^{-1/D})$. If $\ell'$ is outside this range,
the covering has density~$\textbc{k}{t}$, but the probability of this event is
$O(e^{-cl})$, so the total expected number of $k$-sets containing~$T$ coming
from a given hyperplane containing~$T$ is $1+O((\log v)^{-1/D})$, and the total
expected number coming from all such hyperplanes is
$$(1+O((\log v)^{-1/D}))(1+O(v^{-1/t})) = 1+O((\log v)^{-1/D}).$$
Thus the expected density of the induced covering is $1+O((\log v)^{-1/D})$.
\end{proof}

\begin{corollary}
The induced covering algorithm runs in time and space $O(v^t)$.
\end{corollary}

\begin{proof}
By Corollary~\ref{cor:abort}, precomputation takes time $O(\ell^{t+1}\log^2
\ell)$, which is $O(v^t)$ by our choice of~$\ell$. The number of hyperplanes is
$O(p^t) = O(v)$ by our choice of~$p$, so the time to compute the affine
geometry is $O(v^2) = O(v^t)$. For each hyperplane, the work to find the
intersection and convert it into an $(\ell',k,t)$ covering will vary, but the
time per block is constant. Hence the total time and space of the algorithm is
dominated by the size of the $(v,k,t)$ covering, which is also $O(v^t)$.
\end{proof}

\begin{corollary}
The induced covering has expected density
$1 + O \bigl( (\log v)^{-1/D} \bigr)$.
\end{corollary}

Furthermore, if, as we conjecture,
the greedy covering has expected density $1+O(v^{-(k-t)/D})$,
then the expected density of the induced covering  improves to
$1 + O(v^{-(k-t)/D}) + O(v^{-1/t}) = 1 + O(v^{-(k-t)/D})$.

The best way to use the induced covering algorithm in practice is to
first find or make a large table of good coverings with small
parameters using many different methods, and then use these for the
$(\ell',k,t)$ coverings.
We used that strategy to produce
the induced coverings of our earlier paper~\cite{gkp}.

% \bibliography{optimal}
% \bibliographystyle{plain}

\end{document}